\documentclass{ws-m3as}
\usepackage{graphicx}
\begin{document}

\title{ \bf A
MATHEMATICAL MODEL FOR THE EVAPORATION OF \\ A LIQUID FUEL DROPLET,
SUBJECT TO NONLINEAR CONSTRAINTS}

\author{\bf R.Alexandre}
\address{Mathematics Department, University of Evry, 91000 Evry,
France. E-mail: radja.alexandre@univ-evry.fr}


\author{\bf Nguyen Thanh Long}
\address{Department of Mathematics and Computer Science, College of
Natural Science, VietNam National University HoChiMinh City, 227 Nguyen
Van Cu Str., Dist. 5, HoChiMinh City, Vietnam. E-mail: longnt@hcmc.netnam.vn}

\author{\bf A.Pham Ngoc Dinh}
\address{MAPMO, UMR 6628, b\^at. Math\'ematiques, University of Orleans,
BP 6759, 45067 Orl\'eans Cedex 2, France. E-mail:
alain.pham@univ-orleans.fr (corresponding author)}

\maketitle

\begin{abstract} We study the mathematical evolution of a liquid fuel
droplet inside a vessel. In particular, we analyze the evolution of
the droplet radius on a finite time interval. The model problem
involves an hyperbolic system coupled with the pressure and velocity
of the surrounding gas. Existence of bounded solutions for the mass
fraction of the liquid, submitted to nonlinear constraints, is
shown. Numerical simulations are given, in agreement with known
physical experiments.\end{abstract}



\section{Introduction}

\noindent In this paper, we are mainly concerned with the
mathematical analysis of the evaporation of a single droplet in a
gas, in the continuation of our previous work\cite{goutte-1}.\par
\noindent Experimental studies around this subject are of course
important for industrial purposes. Let us refer for instance to the
works\cite{benson,consolini,curtis,hase,jia,lafon,morin,odeide}.\par
\noindent One of our main interest in this paper is to analyze the
time evolution of the droplet radius, a study that we began in our
previous work\cite{goutte-1}. Let us just mention here that the
experimental evolution of this radius is well known as the $d^2$
law\cite{odeide,williams}, where $d$ denotes the diameter of the
droplet, see below for more details.\par
\smallskip
\noindent Let us recall the standard physical framework for this
evolution.\par \noindent The evaporation of a single droplet in a
gas involves simultaneous heat and mass transfer processes. In
particular, heat from evaporation is transferred to the droplet
boundary by conduction and convection, while vapor is carried by
convection and diffusion back into the gas stream. Evaporation rate
depends on the pressure, temperature and physical properties of the
gas, the temperature, volatility and diameter of the drop in the
spray.\par \noindent To fix the ideas, in the experimental
study\cite{odeide} of a single droplet evaporation performed by the
LCSR (Combustion Laboratory of the University of campus from
Orleans, France), the droplet is suspended from a silicate tube. The
elliptic shape of the droplet is assimilated to a sphere of equal
volume. Important quantities of interest for these experiments are
time evolution of the droplet radius, as well as classical
quantities such as mass fractions or temperatures of the liquid and
gas. In the experimental studies performed above, the so-called $
d^{2}$ law is used to simplify two-phase fluid models and then
propose adequate numerical schemes. This law simply states that the
time evolution of the radius behaves as $ {{\displaystyle
d^2}\over{\displaystyle d^{2}_0}}  $ in time flow, and is purely
phenomenological. \par
\smallskip
\noindent Our purpose in this paper is exactly in the opposite
sense. We start from phenomenological fluid (mixtures) PDE modeling
the drop evaporation process, compute the time-evolution of the drop
radius, and then deduce other quantities of interest such as mass
fractions of the liquid and gas. In particular, our numerical
experiments are in good agreement with this phenomenological $d^2$
law, at least for small time evolution.\par
\smallskip
\noindent Our framework is therefore as follows: we consider a
droplet initially represented as a single component mixture (liquid
chemical specie 1) while the surrounding gas at time $ t=0$ is made
of only one (gas) chemical specie, say 2. \par \noindent During the
evaporation process, the liquid vapor is transferred into the gas,
while by condensation at the droplet surface and then by diffusion,
gas chemical specie 2 appears inside the droplet.\ \par \noindent We
make the important simplification that the moving interface between
the droplet and the surrounding gas (i.e. between the two species)
is spherical, with radius $ R=R(t)$ evolving in time.\ \par
\noindent Let $ \rho  _{G}$ (resp. $ \rho  _{L}$) denote gas density
(resp. liquid density), and $ v_{G}$ (resp. $ v_{L}$) denote the gas
velocity (resp. liquid velocity). Then, one has the classical
overall continuity and momentum conservation laws \par
\begin{equation}\label{1.1} \partial  _{t}\rho  _{k}+{\rm div(\rho  }_{k}v_{k}{\rm )}=0
\end{equation}
\begin{equation}\label{1.2}
 \rho  _{k}\partial  _{t}v_{k}+\rho  _{k}v_{k}.\nabla
v_{k}=-\nabla  p .
\end{equation}
Above, subscript $ k$ refers to the gas $ G$ or to the liquid $ L$, depending on whether
one considers the gas or liquid.$ p$ the state equation of the
gas.\ \par
\noindent Let $ Y_{L1},Y_{L2}$ (resp. $ Y_{G1},Y_{G2}$) the mass fractions of
the liquid (resp. gas) obtained after diffusion of species in the surrounding
gas. Therefore for two species, one has\ \par
$$ Y_{G1}+Y_{G2}=Y_{L1}+Y_{L2}=1.$$
\noindent Along with equation (\ref{1.1}), we have to add the equation
giving species conservation. So for the liquid, we have \ \par
\begin{equation}\label{1.3}
 \rho  _{L}\partial  _{t}Y_{Lk}+\rho  _{L}v_{L}.\nabla  Y_{Lk}+{\rm div(\rho
}_{L}Y_{Lk}v_{Lk}{\rm )}=-\rho  _{L}f(Y_{Lk}),\; \; k=1,2 ,
\end{equation}
$ Y_{Lk}$ denoting mass fraction of the liquid, and $ f$ being a continuous
function modeling a friction or a resistance for the drop.\ \par
\noindent We assume that the liquid speed is so small that is can be  settled to $ 0$. Equations (\ref{1.1}) and (\ref{1.3}) can then be written under conservative
form as\ \par
\begin{equation}\label{1.4}
 \partial  _{t}(\rho  \; \tilde g)+{\rm div(\rho  \; }\tilde gv{\rm )}=F(\;
\tilde g)
\end{equation}
or\ \par
\begin{equation}\label{1.5}
 \partial  _{t}u+{{\partial  }\over{\partial  x}}
(f(u))=F(u) ,
\end{equation}
in a system of particular coordinates.\ \par
\noindent If $ \Gamma $ is a curve of discontinuity of $ u$, then one has\ \par
\begin{equation}\label{1.6}
 [f(u)]=[u]{{dx}\over{dt}}  ,
\end{equation}
where $ [.]$ denotes the jump of the inner quantity, $
s={{\displaystyle  dx}\over{\displaystyle
dt}}  $ is the speed of discontinuity along $\Gamma $. The jump relation (\ref{1.6})
is known as Rankine-Hugoniot condition. It merely
means that discontinuities cannot be completely arbitrary. The above considerations are all classical
facts\cite{landau,smoller,ta-tsien,williams}.\ \par
\noindent In the case of our droplet, in order to find interface condition at the droplet surface, i.e. for $ r=R(t)$, it is sufficient to use (\ref{1.5})
and (\ref{1.4}) in polar coordinates, getting\ \par
\begin{equation}\label{1.7}
 [\rho  \; \tilde gv]=[\rho  \; \tilde g]{{dR}\over{dt}}  .
\end{equation}
Thus taking $ \tilde g=1$, one has\ \par
$$ (\rho  _{G}-\rho  _{L}){{dR}\over{dt}}  =\rho  _{G}v_{G}-\rho  _{L}v_{L} ,
$$
that is also with $ v_L=0$\ \par
\begin{equation}\label{1.8}
\rho  _{G}  {\left( v_{G}-{{dR}\over{dt}}  \right) }  =-\rho  _{L}
{{dR}\over{dt}}  .
\end{equation}
Taking $ \tilde g=Y$ in (\ref{1.7}), $ Y$ denoting the mass fraction of the
liquid or the gas after diffusion, we get\ \par
$$ {\left( \rho  _{G}Y_{Gk}-\rho  _{L}Y_{Lk}  \right) }  {{dR}\over{dt}}
=\rho  _{G}Y_{Gk}(v_{G}+v_{Gk})-\rho  _{L}Y_{Lk}v_{Lk}  $$
and this is equivalent to\ \par
\begin{equation}\label{1.9}
 \rho  _{G}Y_{Gk}(v_{G}-R')+\rho  _{G}Y_{Gk}v_{Gk}=-\rho  _{L}Y_{Lk}R'+\rho
_{L}Y_{Lk}v_{Lk}.
\end{equation}
Above $ v_{Gk\  }$(resp.$ v_{Lk}$) is the speed of specie $ Gk$
(resp.$ Lk$), $ k=1,2.$\ \par
\noindent Combining relation (\ref{1.9}) with Fick's law\cite{jia,landau,williams}, that is\ \par
$$ Y_{G1}v_{G1}=_{-}D_{12}\nabla  Y_{G1},\; Y_{G2}v_{G2}=-D_{21}\nabla  Y_{G2} ,
$$
$ D_{12}$ and $ D_{21}$ being diffusion coefficients, and with
equations relating the thermodynamic state
at the interface $ r=R(t)$\ \par
$$ Y_{Gk}=K_{k}Y_{Lk},\; k=1,2 ,$$
we obtain, for the mass fraction of the liquid $ Y_{L1}$,
the boundary condition\ \par
\begin{equation}\label{1.10} \partial  _{r}Y_{L1}+{{R'(t)(K_{1}-1)}\over{K_{2}\rho  _{G}(R(t),t)-K_{3}
}}  Y_{L1}=0\; {\rm at\; }r=R(t) ,
\end{equation}
using polar coordinates.\ \par \noindent In our previous
work\cite{goutte-1}, we have made huge mathematical and physical
simplifications taking the state equation of the gas $ p$ as
constant in (\ref{1.2}), and considering gas velocity $ v_{G}$ as a
given function of the time $ t$. Thus in our previous work, system
(\ref{1.1}), (\ref{1.2}) was reduced to equation (\ref{1.1}) with a
given $ v_{G}(t)$.\ \par \noindent In the present work, we consider
the full hyperbolic system (\ref{1.1}), (\ref{1.2}) with an
auxiliary state equation for the gas given by $ p_{1}=\rho  ^{\gamma
}$. This of course extends our previous work, but considering such
pressure laws has the advantage that we have been able to perform
numerical comparisons. More general state laws will be studied in a
future work.\par \noindent Once $ \rho  _{G}(r,t)$ and $ v_{G}(r,t)$
determined, radius $ R(t)$ of the drop suspended in the gas will be
computed through the ordinary differential equation (\ref{1.8}).
Then we shall determine the mass fraction $ Y_{L1}$ of liquid after
evaporation process, through the PDE (\ref{1.3}) along with boundary
condition (\ref{1.10}), for a given suitable function $ f$.\ \par
\noindent For this last purpose, within the framework of weighted
Sobolev spaces on initial data and for some continuous function $ f$
subject to increasing condition, we shall provide an unique local
solution for the mass fraction $ Y_{L1}$ of the liquid. In addition,
we shall show that if the initial condition is bounded, then so is
our solution.\ \par
 \noindent In the numerical applications (2nd
example) we have chosen the experimental conditions made by the LCSR
in the study of single drop evaporation and in this case the study
of the radius of the drop shows us that the graphic associated to
our mathematical model presents the same features as in the
experimental curves.\ \par
 \noindent{\bf Plan of the paper:} In Section 2,
by using Riemann invariants, we determine the droplet radius. This
enables to get, in Section 3, the liquid mass fraction, using a
variational method. Finally, we have presented some numerical
simulations in the last Section, which shows that our model is is
good agreement with experimental simulations, at least for short
time. \par

\section{Hyperbolic system and droplet radius}

\rm The gas velocity $ v_{G}(r,t)$ and its density $ \rho  _{G}(r,t)$
satisfy the following system, using polar coordinates\ \par
\begin{equation}\label{2.1}
{\matrix {\displaystyle
{{\partial  \rho  _{G}  }\over{\partial  t}}  +{{1}\over{r^{2}  }}  {{\partial
}\over{\partial  r}}  {\left( r^{2}\rho  _{G}v_{G}  \right) }  =0\cr
\displaystyle {{\partial  }\over{\partial  t}}  {\left( \rho
_{G}v_{G}  \right) }  +{{1}\over{r^{2}
}}  {{\partial  }\over{\partial  r}}  {\left( r^{2}\rho  _{G}v_{G}^{2}
\right) }  =-{{\partial  p}\over{\partial  r}} . \cr }}
\end{equation}
Setting $ \rho  (r,t)=r^{2}\rho  _{G}(r,t)$ and $ v(r,t)=v_{G}(r,t)$, we have\ \par
\begin{equation}\label{2.2}
 {\matrix {\displaystyle
{{\partial  \rho  }\over{\partial  t}}  +v{{\partial  \rho  }\over{\partial
r}}  +\rho  {{\partial  v}\over{\partial  r}}  =0\cr
\displaystyle {{\partial  v}\over{\partial  t}}  +v{{\partial
v}\over{\partial  r}}
=-{{1}\over{\rho  }}  {{\partial  p_{1}  }\over{\partial  r}}  ,\cr }}
\end{equation}
where $ p_{1}(r,t)$ is an auxiliary function connected to the state
equation of the gas $ p(r,t)$, by
$$ {{\displaystyle  \partial
p_{1}}\over{\displaystyle  \partial
r}}  =r^{2}{{\displaystyle  \partial  p}\over{\displaystyle  \partial  r}}
.$$
\noindent In (2.2), according to the discussion in the Introduction, we choose the auxiliary function $ p_{1}(r,t)$ as $
p_{1}$=$ \rho  ^{\gamma  },\  \gamma  >1$.\par
\noindent With this choice,
we get
the following system \ \par
\begin{equation}\label{2.3}
 {\matrix {\displaystyle
{{\partial  \rho  }\over{\partial  t}}  +v{{\partial  \rho  }\over{\partial
r}}  +\rho  {{\partial  v}\over{\partial  r}}  =0\cr
\displaystyle {{\partial  v}\over{\partial  t}}  +v{{\partial
v}\over{\partial  r}}
+\gamma  \rho  ^{\gamma  -2}  {{\partial  p}\over{\partial  r}}  =0 .\cr }}
\end{equation}
We note that (\ref{2.3}) is equivalent to matrix form \par
\begin{equation}\label{2.4}
 {{\partial  }\over{\partial  t}}  {\left( {\matrix {
\rho  \cr
v\cr }}  \right) }  +A.{{\partial  }\over{\partial  r}}  {\left( {\matrix {
\rho  \cr
v\cr }}  \right) }  =0,
\end{equation}
$ A$ being the $ (2,2)$  matrix\ \par
$$ A={\left( {\matrix {
v & \rho  \cr
\gamma  \rho  ^{\gamma  -2}   & v\cr }}  \right) } . $$
Eigenvalues of $ A$ (characteristics speeds) are given $ \lambda
=v-c$ and $ \mu  =v+c,\  c={\sqrt{\displaystyle  \gamma  \rho  ^{\gamma  -1}}}
$. Since $ \lambda  <\mu  $, system (\ref{2.4}) is therefore hyperbolic.
Thus there exists two functions $ W(\rho  ,v)$ and $ Z(\rho  ,v)$
(Riemann invariants) such that\ \par
\begin{equation}\label{2.5} W(\rho  ,v)={\rm constant\; on\; }{{dX^{1}  }\over{dt}}  =\lambda  ,
\end{equation}
\begin{equation}\label{2.6}
 Z(\rho  ,v)={\rm constant\; on\; }{{dX^{2}  }\over{dt}}  =\mu
.
\end{equation}
$ W(\rho  ,v)$ is determined by the system $ {{\displaystyle
dv}\over{\displaystyle  c/\rho
}}  ={{\displaystyle d\rho  }\over{\displaystyle 1}}  $, vector
$ R_{1}=(1,c/\rho  )$ being and eigenvector associated to the eigenvalue
$\mu  $. Thus \ \par
\begin{equation}\label{2.7}
 W(\rho  ,v)=v-{{2c}\over{\gamma  -1}}  .
\end{equation}
Similarly, the Riemann invariant $ Z(\rho  ,v)$ corresponding
to $\lambda  $ is given by\ \par
\begin{equation}\label{2.8}
 Z(\rho  ,v)=v+{{2c}\over{\gamma  -1}} .
\end{equation}
Functions $ W(\rho  ,v)=W(t,r)$ and $ Z(\rho  ,v)=Z(t,r)$ satisfy
the following system, equivalent to system (\ref{2.4})\ \par
\begin{equation}\label{2.9}
 {\matrix {\displaystyle
{{\partial  W}\over{\partial  t}}  +\lambda  (W,Z){{\partial  W}\over{\partial
r}}  =0\cr
\displaystyle{{\partial  Z}\over{\partial  t}}  +\mu  (W,Z){{\partial
Z}\over{\partial
r}}  =0 ,\cr }}
\end{equation}
where $ \lambda  (W,Z)$ and $ \mu  (W,Z)$ are given by \ \par
\begin{equation}\label{2.10}
 {\matrix {
\lambda  =-{\left( {{\gamma  -3}\over{4}}  \right) }  Z+{\left( {{\gamma
+1}\over{4}}  \right) }  W\cr
\mu  ={\left( {{\gamma  +1}\over{4}}  \right) }  Z-{\left( {{\gamma
-3}\over{4}}
\right) }  W,\cr }}
\end{equation}
as follows from (\ref{2.7}) and (\ref{2.8}).\ \par
\noindent It is well known that a sufficient condition in order that
(\ref{2.9}) is authentically
nonlinear is that $ {{\displaystyle \partial  \lambda
}\over{\displaystyle \partial
W}}  >0$ and $ {{\displaystyle \partial  \mu  }\over{\displaystyle \partial
Z}}  >0$, which is the case here according to (\ref{2.10}).\ \par
\noindent Integration along the characteristics defined by \ \par
$$ {{dX^{1}  }\over{dt}}  =\lambda  (W,Z),\; X^{1}(0)=\beta  $$
gives \ \par
\begin{equation}\label{2.11}
 X^{1}_{(0,\beta  )}(t)=\beta  +{\int _{0}  ^{t}  {\lambda
(W(s,X^{1}(s)),Z(s,X^{1}(s)))ds}} .
\end{equation}
Therefore, the solution of the initial value problem\ \par
\begin{equation}\label{2.12}
 W_{t}+\lambda  (W,Z)W_{r}=0,\; \; W(0,r)=W_{0}(r)
\end{equation}
can be written as\ \par
\begin{equation}\label{2.13}
 W(t,r)=W_{0}  {\left( X^{1}_{(0,\beta  )}(0)\right) }
=W_{0}(\beta ),
\end{equation}
where $\beta =r-\int_{0}^{t} \lambda (W(s,X^{1}(s)),Z(s,X^{1}(s)))ds$.\ \par
\noindent Similarly we have \ \par
\begin{equation}\label{2.14}
 Z(t,r)=Z_{0}  {\left( X^{2}_{(0,\alpha  )}(0)\right) }  =Z_{0}(\alpha
),
\end{equation}
where\ \par
$$ X^{2}_{(0,\alpha  )}(t)=\alpha  +{\int _{0}  ^{t}  {\mu
(W(s,X^{2}(s)),Z(s,X^{2}(s)))ds}}
.$$
\noindent The above considerations lead to the following
\begin{proposition}\label{prop1} Assume that $ W'_{0}(\beta  )<0\
$or $ Z'_{0}(\alpha  )<0$. Then solution of system (2.3) is
defined on a finite interval $ [0,T[$.\par
\rightline{$\Box$}
\end{proposition}
\noindent Proof: Differentiation of (\ref{2.5}) and (\ref{2.11}) with respect to $\beta  $
gives\ \par
\begin{equation}\label{2.15}
 {{dX^{1}_{\beta  }  }\over{dt}}  =\lambda  _{\beta  }(W,Z)\; \; {\rm with\;
\; }X^{1}_{\beta  }(t=0)=1 ,
\end{equation}
and in the same way\ \par
\begin{equation}\label{2.16} {{dX^{2}_{\alpha  }  }\over{dt}}  =\mu  _{\alpha  }(W,Z)\; \;
{\rm with\; \; }X^{2}_{\alpha  }(t=0)=1.
\end{equation}
Since $ \lambda  _{\beta  }=\lambda  _{W}W_{\beta  }+\lambda
_{Z}Z_{\beta  }=\lambda
_{W}W'_{0}(\beta  )$ and $ \mu  _{\alpha  }=\mu  _{Z}Z'_{0}(\alpha  )$, from (\ref{2.15}) and (\ref{2.16}), integrating w.r.t. $ t$ along
the characteristics yields\ \par
\begin{equation}\label{2.17} X^{1}_{\beta  }(t)=1+{\int _{0}  ^{t}  {\lambda  _{W}W'_{0}(\beta
)dt=1+{\left(
{{\gamma  +1}\over{4}}  \right) }  W'_{0}(\beta  )t}} ,
\end{equation}
\begin{equation}\label{2.18}
 X^{2}_{\alpha  }(t)=1+{\int _{0}  ^{t}  {\mu  _{Z}Z'_{0}(\alpha
)dt=1+{\left(
{{\gamma  +1}\over{4}}  \right) }  Z'_{0}(\alpha  )t}}  .
\end{equation}
\noindent From (\ref{2.17}), it follows that $ X^{1}_{\beta  }(t_{1})=0$ for $
t_{1}={{\displaystyle -4}\over{\displaystyle
(\gamma  +1)W'_{0}(\beta  )}}  >0$. Similarly $ X^{2}_{\alpha  }(t_{2})=0$
for $ t_{2}={{\displaystyle -4}\over{\displaystyle (\gamma  +1)Z'_{0}(\alpha
)}}  >0$.\par
\noindent Hence $ {{\displaystyle \partial  W}\over{\displaystyle \partial
r}}  (t,r)$ becomes infinite for $ T=\displaystyle \inf {\rm
\{}t_{1},t_{2}{\rm \}}$,
since $ {{\displaystyle \partial  W}\over{\displaystyle \partial  r}}
=W_{\beta  }.{{\displaystyle d\beta  }\over{\displaystyle dX}}
={{\displaystyle W'_{0}(\beta
)}\over{\displaystyle X_{\beta  }}}  .$\ \par

\rightline{$\Box$}
\noindent From Proposition \ref{prop1}, it follows that
\begin{proposition}\label{prop-rajout1}
System (\ref{2.9}) admits an unique $ C^{1}$
solution on $ [0,T[$, for all $ r\in  {\mathbb R}^{+}$ and for initial data
$ \rho  _{G}(0,r)=\rho  _{0}(r)$ and $ v_{G}(0,r)=v_{0}(r)$ belonging
to $ C^{1}({\mathbb R}^{+})$.\par
\rightline{$\Box$}
\end{proposition}
\noindent Concerning the droplet radius, it follows from (\ref{1.8}), that we have the following ode for this radius
\begin{equation}\label{2.19}
 {{dR(t)}\over{dt}}  ={{v_{G}(t,R(t))\rho  _{G}(t,R(t))}\over{\rho
_{G}(t,R(t))-\rho
_{L}  }}  ,\; R(0)=R_{0}.
\end{equation}
\noindent We immediately deduce
\begin{proposition}\label{prop-rajout2}
The Cauchy problem (2.19) has an unique solution $ R(t)$ on a maximal
time interval $ [0,T^{*}[$ with $ T^{*}\leq  T$, given initial data
$ \rho  _{0}(r)$ and $ v_{0}(r)$ such that $ W'_{0}(r)<0$ or $
Z'_{0}(r)<0$.\par
\rightline{$\Box$}
\end{proposition}

\section{Liquid Mass Fraction}

\rm The liquid mass fraction $ Y_{L1}$ satisfies the conservation
equation of specie (\ref{1.3}), which can be rewritten as\ \par
\begin{equation}\label{3.1}
 \partial  _{t}Y_{L1}-{{1}\over{r^{2}  }}  {{\partial  }\over{\partial  r}}
{\left( r^{2}  {{\partial  }\over{\partial  r}}  Y_{L1}  \right) }
+f(Y_{L1})=0 .
\end{equation}
We have used polar coordinates, and taken the diffusion constant $ D_{12}$ as being
equal to 1. Of course, (\ref{3.1}) is equivalent to \ \par
\begin{equation}\label{3.2}
 \partial  _{t}Y_{L1}-\Delta Y_{L1}-{{2}\over{r}}  {{\partial
}\over{\partial
r}}  Y_{L1}+f(Y_{L1})=0_{,\; \; }{\rm for\; }0<r<s(t) ,
\end{equation}
where $ s(t)=R(t)$ denotes the droplet radius determined in section 2.\ \par
\noindent The boundary condition at the surface $ s(t)$ is given by
the Rankine-Hugoniot
condition connected to the thermodynamic equilibrium, i.e. formula (\ref{1.10}).\par
\noindent Performing the change of variable $ r=R(t)x$, function $ Y_{L1}(t,r)$
turns to function $ Y_{L1}(t,R(t)x)=u(t,x)$, which satisfies the
following initial boundary value (i.b.v.) problem \ \par
\begin{equation}\label{3.3}
 \partial  _{t}u-a(t){\left( \Delta u+{{2}\over{x}}  \partial  _{x}u\right) }
-x{{R'(t)}\over{R(t)}}  \partial  _{x}u+f(u)=0,\; \; 0<x<1,\;
t>0,
\end{equation}
\begin{equation}\label{3.4}
 {\left\vert {\rm lim}_{x\rightarrow  0+}xu_{x}(t,x)\right\vert }  <\infty
,\; \; u_{x}(t,1)+k(t)u(t,1)=0,
\end{equation}
\begin{equation}\label{3.5}
 u(0,x)=u_{0}(x),
\end{equation}
where we used the following notations
\begin{equation}\label{3.6}
 a(t)={{1}\over{R^{2}(t)}}  ,\; \; k(t)={{R(t)R'(t)(K_{1}-1)}\over{K_{2}\rho
_{G}(t,R(t))-K_{3}  }}  .
\end{equation}
\noindent Our purpose in this Section is to analyze the boundary value problem (\ref{3.3})-(\ref{3.5}).\par
\noindent We shall do so by setting this problem in a variational framework, using weighted Sobolev spaces.\par
\noindent Let $ \Omega =]0,1[$ and define $ H$ as the Hilbert space
given by\ \par
$$ H=\{v:\Omega \rightarrow  {\mathbb R} ,\; {\rm measurable\; and\; such\;
that\; }{\int _{0}  ^{1}  {x^{2}v^{2}(x)dx<+\infty  \}.}}  $$
Note that $ H$ is the closure of $ C^{0}(\bar{\Omega} )$ w.r.t. the
norm $ \Vert
v\Vert _{H}=\Big (\displaystyle \int _{0}^{1}x^{2}v^{2}(x)dx\Big )^{1/2}$.
We also introduce the real Hilbert space $ V=\Big \{v\in  H\  \vert
v'\in  H\Big \}.$
In the following, we shall often use the fact that $ V$ is the closure
of $ C^{1}(\bar{\Omega} )$ w.r.t. the norm $ \Vert v\Vert _{V}=\Big
(\Vert v\Vert
^{2}_{H}+\Vert v'\Vert ^{2}_{H}\Big )^{1/2}$. $ V$ is continuously
embedded in $ H$. Identifying $ H$ with his dual $ H'$, one has $ V\subset
H\subset  V'$with continuous injections.\ \par
\noindent Note also that the norms $ \Vert .\Vert _{H}$ and $ \Vert
.\Vert _{V}$
can be defined, respectively, from the inner products $
<u,v>=\displaystyle \int _{0}^{1}x^{2}u(x)v(x)dx$
and $ <u,v>+<u',v'>$.\ \par
\noindent We then have the following results, the proofs of which can
be found in the paper\cite{bvp},
\begin{lemma}\label{lemma1} \rm For every $ v\in  C^{1}([0,1])$, $
\epsilon  >0$ and
$ x\in  [0,1]$ we have \ \par
$$ \Vert v\Vert ^{2}_{0}\leq{{\displaystyle 1}\over{\displaystyle 2}}
\  \Vert v'\Vert ^{2}_{0}+v^{2}(1) ,$$
$$ v^{2}(1)\leq  \epsilon  \Vert v'\Vert ^{2}_{0}+C_{\epsilon  }\Vert
v\Vert ^{2}_{0} ,$$
$$ \Big \vert v(1)\Big \vert \leq  2\  \Vert v\Vert _{1},\ \Big \vert
xv(x)\Big \vert \leq  {\sqrt{\displaystyle 5}}  \Vert v\Vert _{1}$$
where $ C_{\epsilon  }=3+{{\displaystyle 1}\over{\displaystyle \epsilon  }}
$ and $ \Vert .\Vert _{0}=\Vert .\Vert _{H},\  \Vert .\Vert _{1}=\Vert .\Vert
_{V}$.\ \par
\rightline{$\Box$}
\end{lemma}
\begin{lemma}\label{lemma2} \rm The embedding $ V\subset  H$ is
compact.\ \par
\rightline{$\Box$}
\end{lemma}
\begin{remark}\label{remark1} \rm From Lemma \ref{lemma1}, it follows that $
\Big (\Vert v'\Vert ^{2}_{0}+v^{2}(1)\Big )^{1/2}$
and $ \Vert v\Vert _{1}$ are two equivalent norms on $ V$ since $ {{2}\over{3}}  {\left\Vert v\right\Vert }  ^{2}_{1}\leq
v^{2}(1)+{\left\Vert
v'\right\Vert }  ^{2}_{0}\leq  5{\left\Vert v\right\Vert }  ^{2}_{1}$, for all $ v\in  V$.\par
\rightline{$\Box$}
\end{remark}
\begin{remark} \rm We have $ xv(x)\in  C^{0}([0,1])$, for all $v\in  V$.\par
\noindent Indeed, on one hand, $ \displaystyle \lim _{x\rightarrow  0+}xv(x)=0,\
\forall  v\in  V$ (see the book\cite{adams}, p.128), and on the other hand $ v_{\vert [\epsilon
,1]}\in  C^{0}([\epsilon  ,1]),\  \forall  \epsilon  ,\ 0<\epsilon <1$,
since we have $ H^{1}(\epsilon  ,1)\subset  \  C^{0}([\epsilon  ,1])$
and $ \epsilon  \Vert v\Vert _{H1(\epsilon  ,1)}\leq  \Vert v\Vert _{1\
}\forall  v\in  V$ $ \forall  \epsilon  ,\  0<\epsilon  <1$. \ \par

\rightline{$\Box$}
\end{remark}
\noindent If $X$ is any Banach space, we denote by $ \Vert .\Vert
_{X}$ its norm, and by $ X'$ the dual space of $ X$. We denote by $
L^{p}(0,T;X),\ 1\leq  p\leq  \infty  $, the standard Banach space of
real functions $ u:\ (0,T)\rightarrow  X$, measurable, such that \
\par
$$ {\left\Vert u\right\Vert }  _{  L^{  p}  (0,T;X)}
\; ={\left( {\int _{0}  ^{T}  {{\left\Vert u(t)\right\Vert }  ^{p}_{X}dt}}
\right) }  ^{  1/p}  <+\infty  ,\; \; {\rm for\; }1\leq  p<\infty
$$
and \ \par
$$ {\left\Vert u\right\Vert }  ^{  }  _{  L^{  \infty
}  (0,T;X)}  =\displaystyle{{\rm ess\; sup}_{ 0<t<T}}  {\left\Vert
u(t)\right\Vert }
_{X},\; \; {\rm for\; }p=\infty  .$$
Let $ u(t),\  u'(t),\  u_{x}(t),\  u_{xx}(t)$ denote $ u(t,x),\
{{\displaystyle \partial
u}\over{\displaystyle \partial  t}}  (t,x),\  {{\displaystyle
\partial  u}\over{\displaystyle
\partial  x}}  (t,x),\  {{\displaystyle \partial
^{2}u}\over{\displaystyle \partial
x^{2}}}  (t,x)$ respectively.\ \par
\ \par
\smallskip
\noindent \rm We shall make the following set of \underline{assumptions}:\ \par
\smallskip
\ $\bullet$\ \ (H$_{1}$) $ u_{0}\in  H ;$\ \par
\smallskip
\ $\bullet$\ \ (H$_{2}$) $ a,k\in  W^{1,\infty  }(0,T)$, $ a(t)\geq
a_{0}>0 ;$\ \par
\smallskip
\ $\bullet$\ \ (F$_{1}$) $ f\in  C({\mathbb R}  ,{\mathbb R}  ) ;$\
\par
\smallskip
\ $\bullet$\ \ (F$_{2}$) There exists positive constants $ C_{1},C'_{1},C_{2}$
and $ p,\  1<p<3$, such that\ \par
\ \ \ \ \ \ (i) $ uf(u)\geq  C_{1}\vert u\vert ^{p}-C'_{1} ,$\ \par
\ \ \ \ \ \ (ii) $ \vert f(u)\vert \leq  C_{2}(1+\vert u\vert
^{p-1}) .$\ \par
\medskip
\noindent Let $ u\in  C^{2}([0,T]\times  [0,1])$ be a solution of
problem (\ref{3.3})-(\ref{3.5}).\par
\noindent Then, after multiplying equation (\ref{3.3}) by $ x^{2}v,\  v\in  V$ w.r.t. the
scalar product of $ H$, integrating by parts and taking into
account boundary condition given by (\ref{3.4}), we get\ \par
$$ {{d}\over{dt}}  <u(t),v>+\; a(t){\int _{0}  ^{1}
{x^{2}u_{x}v_{x}dx+a(t)k(t)u(1)v(1)-{{R'(t)}\over{R(t)}}
{\int _{0}  ^{1}  {x^{3}u_{x}vdx+<f(u),v>=0}}  }}  $$
The \underline{weak formulation} of the ibv problem (\ref{3.3})-(\ref{3.5}) can then be given in
the following way: Find $ u(t)$, defined on the open set $ (0,T)$, such that $ u(t)$
satisfies the following variational problem\ \par
\begin{equation}\label{3.7}
{{d}\over{dt}}  <u(t),v>+\; \tilde a(t;u(t),v)+<f(u(t)),v>=0,\; \;
\forall  v\in  V,
\end{equation}
together with the initial condition\ \par
\begin{equation}\label{3.8}
u(0)=u_{0} .
\end{equation}
Above, we have used the following bilinear form
\begin{equation}\label{3.9}
\; \tilde a(t;u,v)=a(t){\int _{0}  ^{1}
{x^{2}u_{x}v_{x}dx+a(t)k(t)u(1)v(1)-{{R'(t)}\over{R(t)}}
{\int _{0}  ^{1}  {x^{3}u_{x}vdx,\; u,v\in  V}}  }}.
\end{equation}
We first note the following lemma, the proof of which can be found in our previous paper\cite{goutte-1}
\begin{lemma}\label{lemma3} \rm There exists constants $ K_{T},\
\alpha  _{T\  }$and
$ \beta  _{T}$ depending on $ T$, such that\ \par
\begin{equation}\label{3.10}
 {\left\vert \; \tilde a(t;u,v)\right\vert }  \leq  K_{T}
{\left\Vert u\right\Vert }
_{1}  {\left\Vert v\right\Vert }  _{1},\; {\rm for\; all\; }u,v\in
V,
\end{equation}
\begin{equation}\label{3.11}
 \; \; \tilde a(t;u,u)\geq  \alpha  _{T}  {\left\Vert u\right\Vert
}  ^{2}_{1}-\beta
_{T}  {\left\Vert u\right\Vert }  ^{2}_{0},\; u,v\in  V.
\end{equation}
\rightline{$\Box$}
\end{lemma}
\noindent We then have the following existence theorem\ \par
\begin{theorem}\label{theorem1}
\rm Let $ T>0$ and assumptions
(H$_{1}$),(H$_{2}$),(F$_{1}$),(F$_{2}$) hold true. Then,
there exists a solution $ u$ of the variational problem (\ref{3.7}),(\ref{3.8}) such that \ \par
$$ u\in  L^{2}(0,T;V)\cap  L^{\infty  }(0,T;H),\  x^{2/p}u\in
L^{p}(Q_{T}) ,$$
$$ tu\in  L^{\infty  }(0,T;V),\  tu_{t}\in  L^{2}(0,T;H).$$
Furthermore, if $ f$ satisfies the additional condition\ \par
$$ (f(u)-f(v))(u-v)\geq  -\delta  \vert u-v\vert ^{2},$$
for all $ u,v\in  {\mathbb R}$, for some $\delta  \in  {\mathbb R}$,
then the above solution $u$ is unique.\ \par
\rightline{$\Box$}
\end{theorem}
\noindent Proof of Theorem \ref{theorem1}. We divide it in several
steps.\ \par
\noindent $\bullet$ Step 1, Galerkin method.\par
\noindent Denote by $ \{w_{j}\},\  j=1,2,...$, an orthonormal basis of the
separable
Hilbert space $ V$. We wish to find $ u_{m}(t)$ of the form \ \par
\begin{equation}\label{3.12}
 u_{m}(t)={\sum _{j=1}  ^{m}  {c_{mj}(t)w_{j}  }} ,
\end{equation}
where $ c_{mj}(t)$ satisfy the following system of nonlinear differential
equations\ \par
\begin{equation}\label{3.13}
<u'_{m}(t),w_{j}>+\; \tilde a(t;u_{m}(t),w_{j})+<f(u_{m}(t)),w_{j}>=0,\;
\; 1\leq  j\leq  m ,
\end{equation}
together with the initial condition\ \par
\begin{equation}\label{3.14}
u_{m}(0)=u_{0m} ,
\end{equation}
and\ \par
\begin{equation}\label{3.15}
u_{0m}\rightarrow  u_{0\; \; \; }{\rm strongly\; in\; }H.
\end{equation}
Clearly, for each $ m$, there exists an unique local solution $ u_{m}(t)$
of the form (\ref{3.12}), which satisfies (\ref{3.13}) and (\ref{3.14}) almost everywhere
on $ 0\leq  t\leq  T_{m}$, for some $ T_{m},\  0<T_{m}\leq  T.$ The
following estimates allow us to take $ T_{m}=T$ for all $ m$.\ \par
\smallskip
\noindent $\bullet$ Step 2, A priori estimates.\ \par
\noindent (a) First estimate. \ \par
\noindent Multiplying $ j^{th}$ equation of system (\ref{3.13}) by
$ c_{mj}(t)$
and summing up w.r.t. $ j$, we have \ \par
\begin{equation}\label{3.16}
 {{1}\over{2}}  {{d}\over{dt}}  {\left\Vert u_{m}(t)\right\Vert }
^{2}_{0}+\;
\; \tilde a(t;u_{m}(t),u_{m}(t))+<f(u_{m}(t)),u_{m}(t)>=0 .
\end{equation}
Using assumption (H$_{2}$), (F$_{2}$,i), Lemma \ref{lemma1} and
Remark \ref{remark1}, it
follows from (\ref{3.16}) that\ \par
\begin{equation}\label{3.17}
 {{d}\over{dt}}  {\left\Vert u_{m}(t)\right\Vert }  ^{2}_{0}+2\alpha  _{T}
{\left\Vert u_{m}(t)\right\Vert }  ^{2}_{1}+2C_{1}  {\int _{0}  ^{1}
{x^{2}  {\left\vert u_{m}(t,x)\right\vert }  ^{p}dx\leq  {{2C'_{1}  }\over{3}}
+2\beta  _{T}  {\left\Vert u_{m}(t)\right\Vert }  ^{2}_{0}  }} .
\end{equation}
Integrating (\ref{3.17}), using (\ref{3.15}), it follows that\ \par
\begin{equation}\label{3.18}
 S_{m}(t)\leq  C_{0}+{{2}\over{3}}  TC'_{1}+2\beta  _{T}  {\int _{0}
^{t}  {S_{m}(s)ds,}}
\end{equation}
where \ \par
\begin{equation}\label{3.19}
 S_{m}(t)={\left\Vert u_{m}(t)\right\Vert }  ^{2}_{0}+2\alpha  _{T}
{\int _{0}  ^{t}  {{\left\Vert u_{m}(s)\right\Vert }
^{2}_{1}ds+2C_{1} {\int _{0}  ^{t}  {ds{\int _{0}  ^{1}  {x^{2}
{\left\vert u_{m}(s,x)\right\vert } ^{p}dx}}  }}  }}  ,
\end{equation}
and $ C_{0}$ is a constant depending only on $ u_{0}$ with $ \Vert u_{0m}\Vert
^{2}_{0}\leq  C_{0\  }\forall  m.$\ \par
\noindent Applying Gronwall's lemma, we obtain from (\ref{3.18})\ \par
\begin{equation}\label{3.20}
 S_{m}(t)\leq  {\left( C_{0}+{{2}\over{3}}  TC'_{1}  \right) }
{\rm exp(}2\beta
_{T}t{\rm )\leq  }M_{T},\; \; \forall  m,\; \forall  t,\; 0\leq  t\leq
T_{m}\leq  T,
\end{equation}
that is $ T_{m}=T$.\ \par
\noindent In the following, we denote by $ M_{T}$  any generic constant depending
only on $ T$.\ \par
\smallskip
\noindent (b) Second estimate.\par
\noindent Replacing $ w_{j}$ by $
t^{2}u_{m}$ in (\ref{3.8}) gives\ \par
\begin{equation}\label{3.21}
 {\matrix {\displaystyle
{\left\Vert tu'_{m}  \right\Vert }^{2}_{0}+{{1}\over{2}}
{{d}\over{dt}} {\left[ a(t){\left\Vert tu_{m}  \right\Vert }
^{2}_{0}+a(t)k(t)t^{2}u^{2}_{m}(1)\right] } +{{1}\over{2}}
{{d}\over{dt}}  {\left[ t^{2}  {\int _{0}  ^{1}  {x^{2}\hat
f(u_{m})dx}}  \right] }  \cr =\displaystyle{\left\Vert u_{mx}
\right\Vert }  ^{2}_{0} {{d}\over{dt}}  {\left[ t^{2}a(t)\right] }
+{{1}\over{2}}  u_{m}^{2}(1){{d}\over{dt}}  (t^{2}a(t)k(t))+2t{\int
_{0} ^{1}  {x^{2}\hat f(u_{m})dx}}  \cr +\displaystyle{{R'(t)t^{2}
}\over{R(t)}}  {\int _{0}  ^{1} {x^{3}u_{mx}u'_{m}dx}} \cr }}
\end{equation}
where \ \par
\begin{equation}\label{3.22}
 \hat f(z)={\int _{0}  ^{z}  {f(y)dy}}  .
\end{equation}
Integrating (\ref{3.21}) w.r.t. time variable from $ 0$ to $ t$, we have,
after some rearrangements\ \par
\begin{equation}\label{3.23}
 {\matrix {\displaystyle
2{\int _{0}  ^{t}  {{\left\Vert su'_{m}(s)\right\Vert }
^{2}_{0}ds+a(t){\left\Vert
tu_{mx}(t)\right\Vert }  ^{2}_{0}+a(t)t^{2}u_{m}^{2}(t,1)}}  \hfill\cr
=\displaystyle a(t)(1-k(t))t^{2}u^{2}_{m}(t,1)+{\int _{0}  ^{t}
{{\left[ s^{2}a(s)\right] }
'{\left\Vert u_{mx}(s)\right\Vert }  ^{2}_{0}ds}} \hfill\cr
+\displaystyle{\int _{0}  ^{t}  {{\left[ s^{2}a(s)k(s)\right] }
'u_{m}^{2}(s,1)ds+2{\int
_{0}  ^{t}  {{{R'(s)}\over{R(s)}}  s^{2}<xu_{mx}(s),u_{m}'(s)>ds}}  }}
\hfill\cr
+\displaystyle 4{\int _{0}  ^{t}  {sds{\int _{0}  ^{1}  {x^{2}\hat
f(u_{m}(s,x))dx-2t^{2}
{\int _{0}  ^{1}  {x^{2}\hat f(u_{m})dx.}} }} }}  \hfill\cr }}
\end{equation}
By means of assumption (H$_{2}$) and Remark \ref{remark1}, we get\ \par
\begin{equation}\label{3.24}
 a(t){\left\Vert tu_{m}(t)\right\Vert }  ^{2}_{0}+a(t)t^{2}u_{m}^{2}(t,1)\geq
{{2}\over{3}}  a_{0}  {\left\Vert tu_{m}(t)\right\Vert }  ^{2}_{1\; \;
\; \; }\forall  t\in  [0,T],\; \forall  m.
\end{equation}
We fix $ \epsilon  >0$ such that \ \par
\begin{equation}\label{3.25}
 {\left\Vert a\right\Vert }  _{\infty  }  {\left\Vert
1-k\right\Vert }  _{\infty
}\epsilon  <{{a_{0}  }\over{3}} ,
\end{equation}
where $ \Vert .\Vert _{\infty  }=\Vert .\Vert _{L\infty  (0,T)}$.\ \par
\noindent Using again Lemma \ref{lemma1}, Remark \ref{remark1} with $
\epsilon  >0$ as in (\ref{3.25}) and
first estimate (\ref{3.20}), the terms on
the r.h.s. of (\ref{3.23}) can be estimated as follows\ \par
\begin{equation}\label{3.26}
 {\matrix {\displaystyle
a(t)(1-k(t))t^{2}u_{m}^{2}(t,1)\leq  {\left\Vert a\right\Vert }  _{\infty  }
{\left\Vert 1-k\right\Vert }  _{\infty  }  {\left( \epsilon
{\left\Vert tu_{m}(t)\right\Vert }
^{2}_{1}+C_{\epsilon  }  {\left\Vert tu_{m}(t)\right\Vert }  ^{2}_{0}
\right) }  \cr
\; \; \; \; \; \; \; \; \; \; \; \; \; \; \; \; \; \; \; \; \; \; \leq
\displaystyle{{a_{0}  }\over{3}}  {\left\Vert tu_{m}(t)\right\Vert }
^{2}_{1\; }+M_{T} ,
\cr }}
\end{equation}
\begin{equation}\label{3.27}
 {\matrix {\displaystyle
{\int _{0}  ^{t}  {{\left( s^{2}a(s)\right) }  '{\left\Vert
u_{mx}(s)\right\Vert }
^{2}_{0}ds+{\int _{0}  ^{t}  {{\left( s^{2}a(s)k(s)\right) }
'u_{m}^{2}(s,1)ds}}
}}  \hfill\cr
\displaystyle \leq  {\left[ {\left\Vert (t^{2}a)'\right\Vert }
_{\infty  }+{\left\Vert (t^{2}ak)'\right\Vert }
_{\infty  }  \right] }  {\int _{0}  ^{t}  {{\left[ {\left\Vert
u_{mx}(s)\right\Vert }
^{2}_{0}+u_{m}^{2}(s,1)\right] }  ds}} \hfill \cr
\displaystyle \leq  5{\left[ {\left\Vert (t^{2}a)'\right\Vert }
_{\infty  }+{\left\Vert (t^{2}ak)'\right\Vert }
_{\infty  }  \right] }  {\int _{0}  ^{t}  {{\left\Vert
u_{m}(s)\right\Vert }  ^{2}_{1}ds\leq
M_{T}  }} , \cr }}
\end{equation}
\begin{equation}\label{3.28}
 2{\left\vert {\int _{0}  ^{t}  {s^{2}  {{R'(t)}\over{R(t)}}
<xu_{mx}(s),u'_{m}(s)>ds}}
\right\vert }  \leq  {\int _{0}  ^{t}  {{\left\Vert
su_{m}'(s)\right\Vert }  ^{2}_{0}ds+{\left\Vert
{{R'}\over{R}}  \right\Vert }  ^{  2}  _{  \infty  }
{\int _{0}  ^{t}  {{\left\Vert su_{m}(s)\right\Vert }  ^{2}_{1}ds.}}  }}
\end{equation}
From assumptions (F$_{1}$) and (F$_{2}$), we note also that\ \par
\begin{equation}\label{3.29}
 -\hat m_{0}=-{\int _{-z_{0}  }  ^{z_{0}  }  {{\left\vert
f(y)\right\vert }  dy\leq
\hat f(z)={\int _{0}  ^{z}  {f(y)dy\leq  C_{2}  {\left( {\left\vert
z\right\vert }
+{{{\left\vert z\right\vert }  ^{p}  }\over{p}}  \right) }  ,\; \; \forall
z\in  {\mathbb R}  }}  }} ,
\end{equation}
where $ z_{0}=(C_{1}'/C_{1})^{1/p}.$\ \par
\noindent Using first estimate (\ref{3.20}), (\ref{3.29}) and Lemma
\ref{lemma1}, we obtain \ \par
\begin{equation}\label{3.30}
 {\matrix {\displaystyle
{\left\vert 4{\int _{0}  ^{t}  {sds{\int _{0}  ^{1}  {x^{2}\hat
f(u_{m}(s,x))dx-2t^{2}
{\int _{0}  ^{1}  {x^{2}\hat f(u_{m}(t,x)dx}}  }}  }}  \right\vert }\hfill \cr
\displaystyle\leq  4C_{2}  {\int _{0}  ^{t}  {sds{\int _{0}  ^{1}
{x^{2}  {\left( {\left\vert
u_{m}(s,x)\right\vert }  +{{1}\over{p}}  {\left\vert
u_{m}(s,x)\right\vert }  ^{p}
\right) }  dx+2t^{2}  {\int _{0}  ^{1}  {x^{2}\hat m_{0}dx}}  }}  }}
\hfill\cr
\displaystyle\leq  4C_{2}  {\int _{0}  ^{t}  {s{\left\Vert
u_{m}(s)\right\Vert }  _{0}ds+{{4}\over{p}}
C_{2}t{\int _{0}  ^{t}  {ds{\int _{0}  ^{1}  {x^{2}  {\left\vert
u_{m}(s,x)\right\vert }
^{p}dx+{{2}\over{3}}  T^{2}\hat m_{0}  }}  }}  }}  \hfill\cr
\displaystyle\leq  2C_{2}T{\sqrt{ M_{T}  }}  +{{2C_{2}  }\over{pC_{1}
}}  TM_{T}+{{2}\over{3}}
T^{2}\hat m_{0}\leq  M_{T}.  \cr }}
\end{equation}
Hence, we deduce from (\ref{3.23}), (\ref{3.24}), (\ref{3.26})-(\ref{3.28}) and (\ref{3.30}) that\ \par
\begin{equation}\label{3.31}
 {\int _{0}  ^{t}  {{\left\Vert su'_{m}(s)\right\Vert }  ^{2}_{0}ds+{{a_{0}
}\over{3}}  {\left\Vert tu_{m}(t)\right\Vert }  ^{2}_{1}\leq  M_{T}+{\left\Vert
{{R'}\over{R}}  \right\Vert }  ^{  2}  _{  \infty  }
{\int _{0}  ^{t}  {{\left\Vert su_{m}(s)\right\Vert }  ^{2}_{1}ds}}  }}
.
\end{equation}
By Gronwall's lemma, we get \ \par
\begin{equation}\label{3.32}
{\int _{0}  ^{t}  {{\left\Vert su'_{m}(s)\right\Vert }  ^{2}_{0}ds+{{a_{0}
}\over{3}}  {\left\Vert tu_{m}(t)\right\Vert }  ^{2}_{1}\leq
M_{T}{\rm exp}{\left(
{\left\Vert {{R'}\over{R}}  \right\Vert }  ^{  2}  _{  \infty
}  {\rm .}{{3T}\over{a_{0}  }}  \right) }  \leq  M_{T},\; \; \; \forall  t\in
[0,T].}}
\end{equation}
Finally, using (\ref{3.20}) and assumption (F$_{2}$,ii) we have also\ \par
\begin{equation}\label{3.33}
 {\int _{0}  ^{t}  {ds{\int _{0}  ^{1}  {{\left\vert x^{  2/p'}
f(u_{m}(s,x))\right\vert }  ^{  p'}  dx\leq  (2C_{2})^{  p'}
{\int _{0}  ^{t}  {ds{\int _{0}  ^{1}  {x^{2}  {\left\vert
u_{m}(s,x)\right\vert }
^{  p}  dx\leq  M_{T}  }}  }}  }}  }} ,
\end{equation}
with $ p'={{\displaystyle p}\over{\displaystyle p-1}}  $.\ \par
\smallskip
\noindent $\bullet$ Step 3, the limiting process.\par
\noindent From (\ref{3.20}), (\ref{3.32}) and (\ref{3.33}), we deduce that there exists a
subsequence
of $ \{u_{m}\}$, still denoted $ \{u_{m}\}$ such that\ \par
\begin{equation}\label{3.34}
 {\matrix {
u_{m}\rightarrow  u\; {\rm weakly\; }*\; {\rm in\; }L^{\infty  }(0,T;H)\; ,
\hfill\cr
u_{m}\rightarrow  u\; {\rm weakly\; in\; }L^{2}(0,T;V)\; ,\hfill\cr
x^{2/p}u_{m}\rightarrow  x^{2/p}u\; {\rm weakly\; in\; }L^{p}(Q_{T})\;
^{  }  \; ,\hfill\cr
tu_{m}\rightarrow  tu\; {\rm weakly\; }*\; {\rm in\; }L^{\infty  }(0,T;V)\; ,
\hfill\cr
(tu_{m})'\rightarrow  (tu)'{\rm weakly\; \; in\; }L^{2}(0,T;H).\;
\; \hfill\cr
\; \cr }}
\end{equation}
Using a standard compactness lemma\cite{lions} (p.57) together with (\ref{3.34}), we can extract from the sequence $ \{u_{m}\}$, a subsequence still denoted by $ \{u_{m}\}$
such that\ \par
\begin{equation}\label{3.35}
 tu_{m}\rightarrow  tu\; \; {\rm strongly\; in\;
}L^{2}(0,T;H).
\end{equation}
Continuity of $ f$ also implies (up to a sub-sequence)\ \par
\begin{equation}\label{3.36}
 f(u_{m}(t,x))\rightarrow  f(u(t,x))\; \; \; {\rm a.e.\; }(t,x)\in
Q_{T}=(0,T)\times  (0,1) .
\end{equation}
\noindent Applying a standard weak convergence lemma\cite{lions}, we have also\ \par
$$ x^{  2/p'}  f(u_{m})\rightarrow  x^{  2/p'}
f(u)\; {\rm weakly\; in\; }L^{  p'}  (Q_{T}).$$
Passing to the limit in (\ref{3.13}) and (\ref{3.14}), it follows from (\ref{3.15}), (\ref{3.34}) and
(\ref{3.36}), that function $ u(t)$ satisfies the i.b.v. problem
(\ref{3.7}), (\ref{3.8}).\ \par
\smallskip
\noindent $\bullet$ Step 4. Uniqueness of the
solutions.\ \par
\noindent First of all, we note the following slight extension of a lemma used in our previous paper\cite{goutte-1} (see also the book\cite{lions})
\begin{lemma}\label{lemma4}\rm Let $ w$ be the weak solution of the
following i.b.v. problem\ \par
\ \ \ $ w_{t}-a(t)(w_{xx}+{{\displaystyle 2}\over{\displaystyle x}}  w_{x})=\
\tilde f(t,x),\  0<t<T,\  0<x<1,$\ \par
\ \ \ $ \Big \vert \displaystyle \lim _{x\rightarrow  0+}xw_{x}(t,x)\Big \vert
<+\infty  ,\  \  w_{x}(t,1)+k(t)w(t,1)=0,\  w(0,x)=0,$ \ \par
\ \ \ $ w\in  L^{2}(0,T;V)\cap  L^{\infty  }(0,T;H)$\it , $ x^{2/p}w\in
L^{p}(Q_{T}) ,$\ \par
\ \ \ $ tw\in  L^{\infty  }(0,T;V)$, $ tw_{t}\in  L^{2}(0,T;H) .$\ \par
\noindent \rm Then\ \par
$$ {\matrix {\displaystyle
{{1}\over{2}}  {\left\Vert w(t)\right\Vert }  ^{2}_{0}+{\int _{0}  ^{t}
{a(s){\left[ {\left\Vert w_{x}(s)\right\Vert }
^{2}_{0}+k(s)w^{2}(s,1)\right] }
ds-{\int _{0}  ^{t}  {<\tilde f(s),w(s)>ds=0,\; }}  }}  \cr
\; \; \; \; \; \; \; \; \; \; \; \; \; \; \; \; \; \; \; \; \; \; \;
\; \; \; \; \; \; \; \; \; \; \; \; \; \; \; \; \; \; \; \; \; \; \;
{\rm a.e.\; }t\in  (0,T)\cr }}  $$
\rightline{$\Box$}
\end{lemma}
\noindent Uniqueness of solutions for our initial i.b.v problem will then be deduced as follows. Let
$ u$ and $ v$
be two weak solutions of (\ref{3.3})-(\ref{3.5}). Then $ w=u-v$ is a weak solution
of problem mentioned in Lemma \ref{lemma4}, with r.h.s. given by $ \tilde f(t,x)={{\displaystyle
xR'(t)}\over{\displaystyle R(t)}}  w_{x}-f(u)+f(v)$. Therefore, Lemma
\ref{lemma4} implies \ \par
$$ {{1}\over{2}}  {\left\Vert w(t\right\Vert }  ^{2}_{0}+{\int _{0}
^{t}  {\; \; \tilde a(s;w(s),w(s))ds+2{\int _{0}  ^{t}
{<f(u(s))-f(v(s)),w(s)>ds=0^{
}  }}  }}  .$$
Using Lemma \ref{lemma3} and assumption (F$_{3}$) we obtain\ \par
\begin{equation}\label{3.37}
 {\left\Vert w(t)\right\Vert }  ^{2}_{0}+2\alpha  _{T}  {\int _{0}
^{t}  {{\left\Vert w(s)\right\Vert }  ^{2}_{1}\leq  2(\delta  +\beta
_{T}){\int
_{0}  ^{t}  {{\left\Vert w(s)\right\Vert }  ^{2}_{0}ds}}  }}
.
\end{equation}
If $ \delta  +\beta  _{T}\geq  0$ we have $ \Vert w(t)\Vert _{0}=0$
by applying Gronwall's lemma. In the case where $ \delta  +\beta  _{T}<0$
the result is clearly still true. \ \par
\noindent This ends the proof of Theorem \ref{theorem1}.\par
\rightline{$\Box$}
\smallskip
\noindent We now turn to the boundness of the above solutions.\par
\noindent For this purpose, we shall make use of the following
assumptions\ \par
\smallskip
\ $\bullet$\ \ (H'$_{1}$) $ u_{0}\in  L^{\infty  }(0,1)$, $ \vert
u_{0}(x)\vert \leq
M,\  a.e.\  x\in  (0,1)$\ \par
\smallskip
\ $\bullet$\ \ (H'$_{2}$) $ a,k\in  W^{1,\infty  }(0,T)$, $ a(t)\geq  a_{0}>0$,
$ k(t)\geq  k_{0}>0$\ \par
\smallskip
\ $\bullet$\ \ (F'$_{1}$) $ uf(u)\geq  0\  \  \forall  u\in  {\mathbb R}
   \  {\rm such\
that\  }\vert u\vert \geq  \Vert u_{0}\Vert _{\infty  }$, for a.e.
$ x\in  (0,1).$\ \par
\smallskip
\noindent We then have the following result\ \par
\begin{theorem}\label{theorem2}
\rm Let (H'$_{1}$), (H'$_{2}$), (F$_{1}$)-(F$_{3}$) and (F'$_{1}$)
hold. Then the unique weak solution of the ibv problem (\ref{3.7})-(\ref{3.9}),
as given by theorem 1, belongs to $ L^{\infty  }(Q_{T})$.\ \par

\rightline{$\Box$}
\end{theorem}
\noindent Proof of Theorem \ref{theorem2}. Firstly, we note that $Z=u-M$ satisfies the i.b.v. problem\ \par
\begin{equation}\label{3.38}
 \partial  _{t}Z-a(t){\left( \Delta Z+{{2}\over{x}}  \partial  _{x}Z\right) }
-x{{R'(t)}\over{R(t)}}  \partial  _{x}Z+f(Z+M)=0,\; 0<x<1,\; t\in  (0,T) ,\;
\end{equation}
\begin{equation}\label{3.39}
 {\left\vert {\rm lim}_{x\rightarrow  0+}xZ_{x}(t,x)\right\vert }  <\infty
,\; \; Z_{x}(t,1)+k(t){\left[ Z(t,1)+M\right] }  =0 ,
\end{equation}
\begin{equation}\label{3.40}
 Z(0,x)=u_{0}(x)-M
\end{equation}

\noindent Multiplying equation (\ref{3.38}) by $ x^{2}v$, for $ v\in  V$, integrating
by parts w.r.t. variable $ x$ and taking into account boundary condition
(\ref{3.39}), one has\ \par
\begin{equation}\label{3.41}
 {\matrix {\displaystyle
{\int _{0}  ^{1}  {x^{2}Z_{t}vdx+a(t){\int _{0}  ^{1}
{x^{2}Z_{x}v_{x}dx+a(t)k(t)Z(t,1)v(1)-{{R'(t)}\over{R(t)}}
{\int _{0}  ^{1}  {x^{3}Z_{x}vdx}}  }}  }}  \cr
\displaystyle +{\int _{0}  ^{1}  {x^{2}f(Z+M)vdx=-Ma(t)k(t)v(1),\; \;
\forall  v\in
V}} , \cr }}
\end{equation}
hence for $ v=Z^{+}={{\displaystyle
1}\over{\displaystyle 2}}
\Big (Z+\vert Z\vert \Big )$, since $ u_{0}\in  L^{\infty  }(0,1)$.
It follows that\ \par
\ \ \ \ \ $ {{\displaystyle 1}\over{\displaystyle 2}}
{{\displaystyle d}\over{\displaystyle
dt}}  \displaystyle \int _{0}^{1}x^{2}\vert Z^{+}\vert
^{2}dx+a(t)\displaystyle \int
_{0}^{1}x^{2}\vert (Z^{+})_{x}\vert ^{2}dx+a(t)k(t)\vert Z^{+}(t,1)\vert
^{2}$\ \par
\ \ \ \ \ \ $ -{{\displaystyle R'(t)}\over{\displaystyle R(t)}}
\displaystyle \int _{0}^{1}x^{3}Z^{+}_{x}Z^{+}dx+\displaystyle \int
_{0}^{1}x^{2}f(Z^{+}+M)Z^{+}dx=-Ma(t)k(t)Z^{+}(t,1)\leq  0 ,$\ \par
\noindent since \ \par
\ \ \ \ \ \ $ \displaystyle \int
_{0}^{1}x^{2}Z_{t}Z^{+}dx=\displaystyle \int
_{0,Z>0}^{1}x^{2}(Z^{+})_{t}Z^{+}dx={{\displaystyle
1}\over{\displaystyle 2}}  {{\displaystyle d}\over{\displaystyle dt}}
\displaystyle \int _{0}^{1}x^{2}\vert Z^{+}\vert ^{2}dx$.\ \par
\noindent On the other hand, by assumption (H'$_{2}$) and Remark
\ref{remark1}, one has\ \par
\begin{equation}\label{3.42}
 a(t){\int _{0}  ^{1}  {x^{2}  {\left\vert Z^{+}_{x}  \right\vert }
^{2}dx+a(t)k(t){\left\vert
Z^{+}(t,1)\right\vert }  ^{2}\geq  \; \tilde C_{0}  {\left\Vert
Z^{+}(t)\right\Vert }
^{2}_{1},}}
\end{equation}
where $ \tilde C_{0}={{\displaystyle 2}\over{\displaystyle 3}}  a_{0}\min
{\rm \{}1,k_{0}{\rm \}}$. \ \par
\noindent Using the monotonicity of $ f(u)+\delta  u$ and (F'$_{1}$),
we have\ \par
\begin{equation}\label{3.43}
 {\matrix {\displaystyle
{\int _{0}  ^{1}  {x^{2}f(Z^{+}+M)Z^{+}dx={\int _{0}  ^{1}  {x^{2}  {\left[
f(Z^{+}+M)-f(M)\right] }  Z^{+}dx+{\int _{0}  ^{1}  {f(M)x^{2}Z^{+}dx}}
}}  }}  \cr
\; \; \; \; \; \; \; \; \displaystyle\geq  -\delta  {\int _{0}  ^{1}
{x^{2}  {\left\vert Z^{+}
\right\vert }  ^{2}dx+{\int _{0}  ^{1}  {f(M)x^{2}Z^{+}dx\; \geq  -\delta
{\left\Vert Z^{+}  \right\Vert }  ^{2}_{0}.}}  }}  \cr }}
\end{equation}
(\ref{3.41})-(\ref{3.43}) together with Cauchy's inequality applied
to the term $ -{{\displaystyle R'(t)}\over{\displaystyle R(t)}}
\displaystyle \int _{0}^{1}x^{3}Z^{+}_{x}Z^{+}dx$ yields\ \par
\begin{equation}\label{3.44}
 {{d}\over{dt}}  {\left\Vert Z^{+}(t)\right\Vert }  ^{2}_{0}+\tilde C_{0}
{\left\Vert Z^{+}(t)\right\Vert }  ^{2}_{1\; }\leq  \; {\left( {{1}\over{\tilde
C_{0}  }}  {\left\Vert {{R'}\over{R}}  \right\Vert }  ^{2}_{\infty
}+2{\left\vert
\delta  \right\vert }  \right) }  {\left\Vert Z^{+}(t)\right\Vert }
^{2}_{0}.
\end{equation}
Integrating (\ref{3.44}), we get\ \par
\begin{equation}\label{3.45}
 {\left\Vert Z^{+}(t)\right\Vert }  ^{2}_{0}\leq  {\left\Vert
Z^{+}(0)\right\Vert }
^{2}_{0}+\; {\left( {{1}\over{\tilde C_{0}  }}  {\left\Vert {{R'}\over{R}}
\right\Vert }  ^{2}_{\infty  }+2{\left\vert \delta  \right\vert }  \right) }
{\int _{0}  ^{t}  {{\left\Vert Z^{+}(s)\right\Vert }  ^{2}_{0}ds.}}
\end{equation}
Since $ Z^{+}(0)=\Big (u(0,x)-M\Big )^{+}=\Big (u_{0}(x)-M\Big )^{+}=0$, Gronwall's lemma yields $ \Vert Z^{+}(t)\Vert _{0}=0.$
Thus $ u(t,x)\leq  M\  $a.e. $ (t,x)\in  Q_{T}.$ \ \par
\noindent The case $ u_{0}(x)\geq  -M$ is similar, by considering $ Z=u+M$ and $ Z^{-}={{\displaystyle
1}\over{\displaystyle 2}}
\Big (\vert Z\vert -Z\Big )$. Thus we get $ Z^{-}=0$ and hence $ u(t,x)\geq
-M$ a.e. $ (t,x)\in  Q_{T}$. \ \par
\noindent All in all, one obtains $ \vert u(t,x)\vert \leq  M\  $a.e.
$ (t,x)\in
Q_{T}$ and this ends the proof of Theorem \ref{theorem2}.\ \par
\rightline{$\Box$}

\section{Numerical applications}

\noindent \rm For the numerical applications, we have taken in (\ref{2.3})
$ \gamma  =3$,
so that equation (\ref{2.9}) reduces to Burger's equation\ \par
\begin{equation}\label{4.1}
 {\left\lbrace {\matrix {
W_{t}+WW_{r}=0,\; \; W(0,r)=W_{0}(r)\cr
Z_{t}+ZZ_{r}=0,\; \; Z(0,r)=Z_{0}(r).\cr }}  \right.}
\end{equation}
It is well known that classical Burger's equation
$$ u_{t}+uu_{r}=0,\; \; u(0,r)=u_{0}(r)$$
admits the solution $ u(t,r)=u_{0}(\xi  (t,r))$, $ \xi  (t,r)$
being defined by the parametrization $ r=u_{0}(\xi  )t+\xi$. \ \par
\smallskip
\noindent Having in mind (\ref{4.1}), we have considered two examples.\par
\smallskip
\noindent $\bullet$ First example.\par
\noindent For the first
example, we have chosen the initial conditions $ W_{0}(r)=1,\  r>0;\
W_{0}(0)=0$ and $ Z_{0}(r)=2,$ $ r>0;$ $ Z_{0}(0)=0$. The continuous
solutions of (\ref{4.1}) are then given by\ \par
$$ W(t,r)={\left\lbrace {\matrix {
1\; {\rm if\; }0\leq  t\leq  r ,\cr
{{r}\over{t}}  \; {\rm if\; }0\leq  r\leq  t ,\cr }}  \right.}  $$
$$ Z(t,r)={\left\lbrace {\matrix {
2\; {\rm if\; }0\leq  2t\leq  r ,\cr {{r}\over{t}}  \; {\rm if\;
}0\leq  r\leq  2t.\cr }}  \right.}  $$ According to Section 2, the
droplet radius is given by formula (\ref{2.19}) which we consider
here with an initial condition taken equal to be $ 1$
\begin{equation}\label{4.2}
 {{dR(t)}\over{dt}}  ={{v_{G}(t,R(t))\rho  _{G}(t,R(t))}\over{\rho
_{G}(t,R(t))-\rho
_{L}  }}  ,\; R(0)=1 ,
\end{equation}
and where
\begin{equation}\label{4.3}
 {\left\lbrace {\matrix {\displaystyle
v_{G}(t,r)={{1}\over{2}}  {\left( W(t,r)+Z(t,r)\right) }  \cr
\displaystyle\rho  _{G}(t,r)={{1}\over{2{\sqrt{ 3}}  r^{2}  }}
{\left( Z(t,r)-W(t,r)\right) }
.\cr }}  \right.}
\end{equation}
In figure 1 below, we have drawn the curve $ t\longrightarrow  R(t)$ on the
time interval $ [0,1]$ with a step $ h=0.05$ and $ \rho  _{L}=0.9.$\ \par
\par
\centerline{
\includegraphics[width=3in]{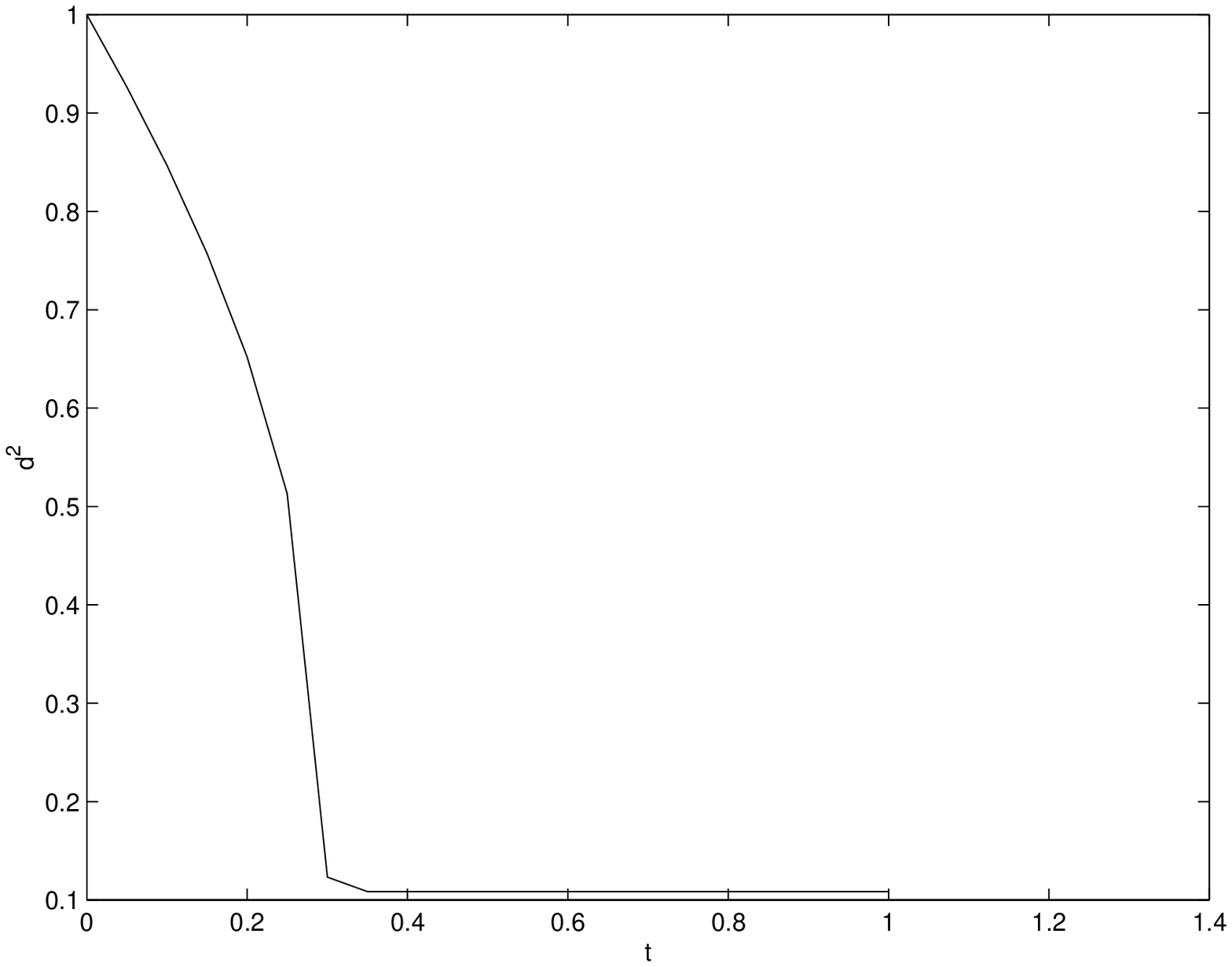}}
\centerline{Fig. 1}\par
\smallskip
\noindent $\bullet$ Second example.\par \noindent  For the second
example, we have chosen the truly experimental conditions made by
the LCSR in the study of single drop evaporation, the drop being
suspended from a silicate tube. Drops are made up of n-heptane fuel
($ \rho  _{L}=683\  kg/mm^{3})$) in air at normalized atmospheric
pressure and with an initial speed $ v_{G}(0,r)=C_{1}=35mm/s$. The
initial density $ \rho  _{G}(0,r)$ of the gas is taken as $
C_{2}={{\displaystyle 348}\over{\displaystyle T_{0}}}  ,\
T_{0}=373K$.\par \noindent In this case, the solution of (\ref{4.1})
are given by \ \par \ \ \ $ W(t,r)=C_{1}-{\sqrt{\displaystyle 3}}
C_{2}\xi  ^{2},\  \xi ={{\displaystyle 1+{\sqrt{\displaystyle
1-4(r-C_{1}t){\sqrt{\displaystyle 3}}  C_{2}t}} }\over{\displaystyle
2{\sqrt{\displaystyle 3}}  C_{2}t}}  $,\ \par \ \ \ $
Z(t,r)=C_{1}+{\sqrt{\displaystyle 3}}  C_{2}\eta  ^{2},\  \eta
={{\displaystyle -1+{\sqrt{\displaystyle
1+4(r-C_{1}t){\sqrt{\displaystyle 3}} C_{2}t}}  }\over{\displaystyle
2{\sqrt{\displaystyle 3}}  C_{2}t}}  $. \ \par \noindent We then
compute $v_G$ and $\rho_G$ by formula (\ref{4.3}), and then solve
the ode for $R$ given by (\ref{4.2}).\par \noindent We note that
$R'(0) = {{C_1C_2}\over{C_2 -\rho_L}} <0 \mbox{ where } C_2 =
{{348}\over{372}} < 1 < \rho_L = 683$. \par
 \noindent On the other hand, one has
$$v_G (t,r) = {1\over{4\sqrt 3 C_2 t}} [ C_1  + {{\sqrt 3}\over{12 t^2 C_2}} (\alpha + \sqrt{1-\alpha} -\sqrt{1+\alpha}) ] ,$$
where
$$\alpha =\alpha (t,r) = 4 (r-C_1 t )\sqrt 3 C_2 t , \ | \alpha | \leq 1 .$$
One has $\mbox{ sgn} v_G (t,r) = \mbox{ sgn} [ C_1 12 t^2 C_2 +\sqrt
3 (\alpha +\sqrt{ 1 -\alpha} -\sqrt{1+\alpha})]$. Since $|\alpha |
\leq 1$, it follows that $\sqrt 3 (\alpha + \sqrt{1-\alpha}
-\sqrt{1+\alpha} ) \geq -\sqrt 6$.\par
 \noindent Thus  $v_G (t,r)>0$ if $C_1 12 t^2 C_2 \geq \sqrt 6 $ i.e. $ t \geq t_2 = {1\over{C_1 2^{1\over 4}}} \simeq
 {1\over{35}}$ and $v_G (t,r)<0$ on the interval $[t_1,t_2[$ .\par
 \noindent Similarly $0\leq\rho_G (t,r) \leq {1\over{C^2_2 t^2 6 .2.r^2}} [ 2 +\sqrt{1-\alpha}
-\sqrt{1+\alpha}] \leq {1\over{3t^2r^2}}$.\par
 \noindent Thus, if $t\geq t_m$, $r\geq r_m$, $\rho_G (t,r) \leq {1\over{3t^2_m r^2_m}}
\leq \rho_L \simeq 683$, that is for $t_m r_m \geq
{1\over{3.683}}\simeq {1\over{45}}$.\par \noindent Since $R'(0)<0$
it follows that $R(t)$ is decreasing on $[0,t_1[$, increasing on
$[t_1,t_2[$ and then from the starting point $t_2$ always non
increasing.\par \noindent Figures 2 and 3 represent resp. the
velocity $ v_{G}(t,r)$ and the pressure $ \rho _{G}(t,r)$ given by
(4.3) for $ (t,r)\in (0,1)\times (0,1)$. \
\par
\par
\centerline{
\includegraphics[width=2.75in]{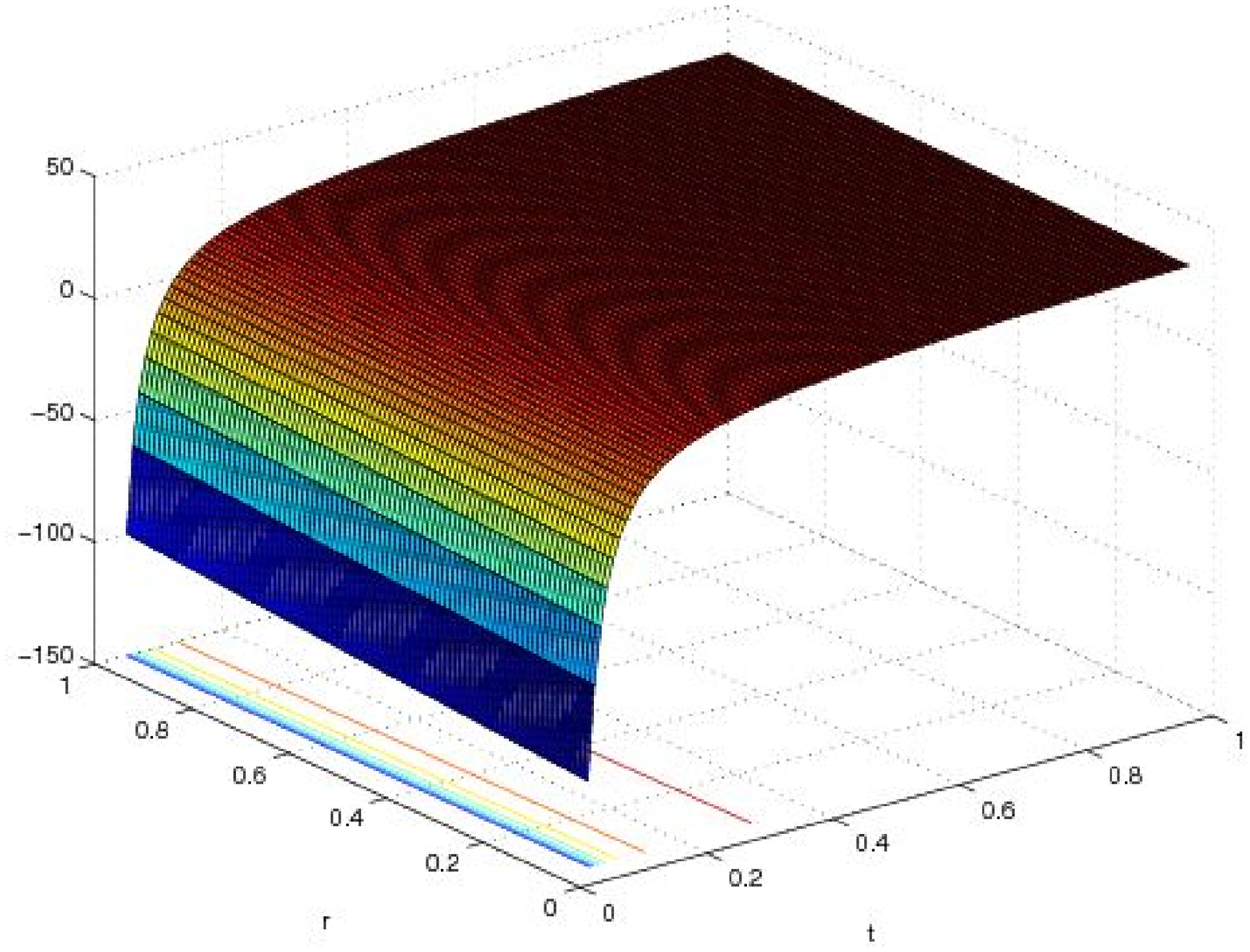}}
\centerline{Fig. 2}\par  \centerline{ \vspace {0.5cm}
\includegraphics[width=2.75in]{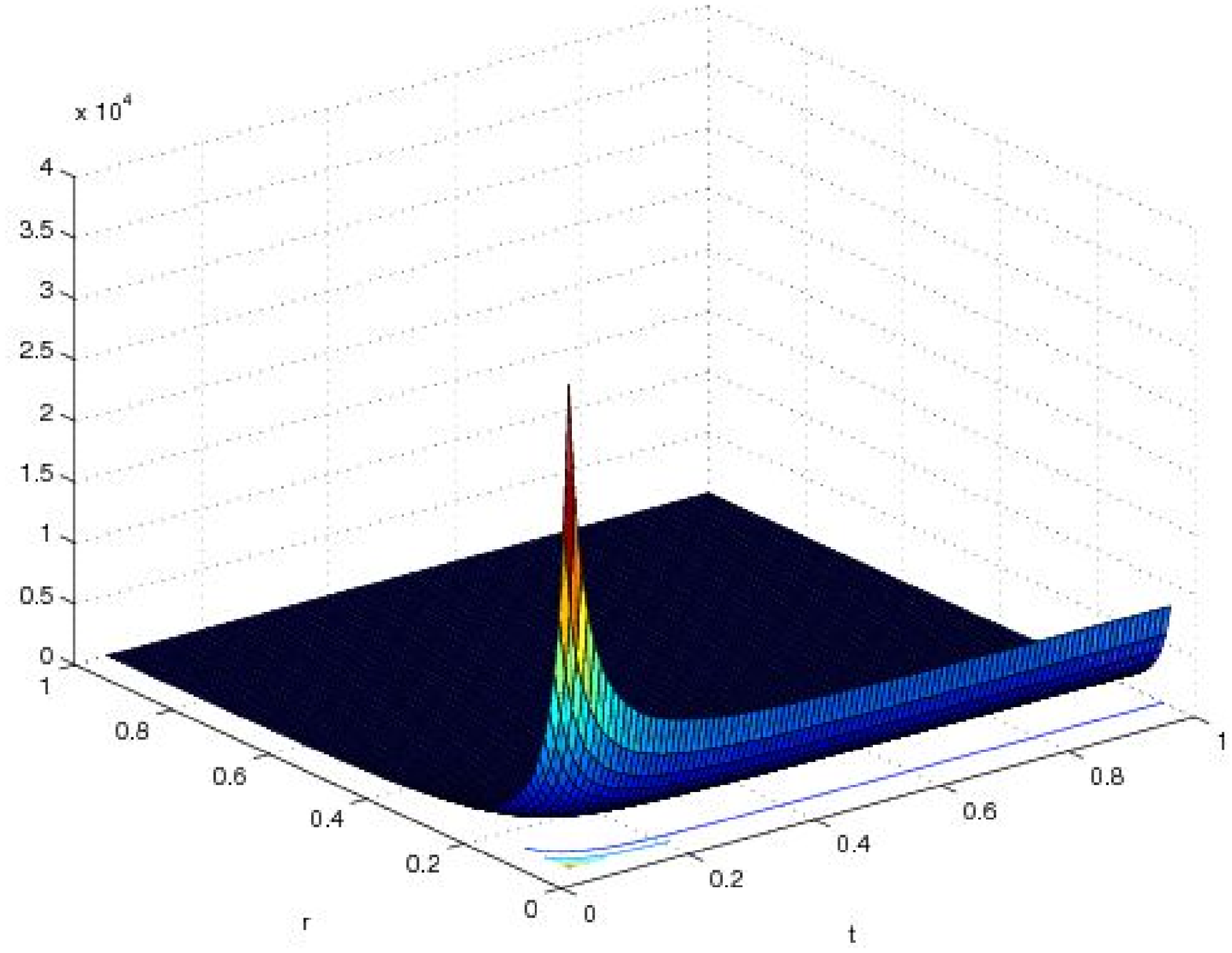}}
\centerline{Fig. 3}\par \vspace {1cm}The curve of the radius $
t\longrightarrow R(t)$ for this case is drawn in figure 4. \ \par
\vspace {0.5cm}
\centerline{
\includegraphics[width=3in]{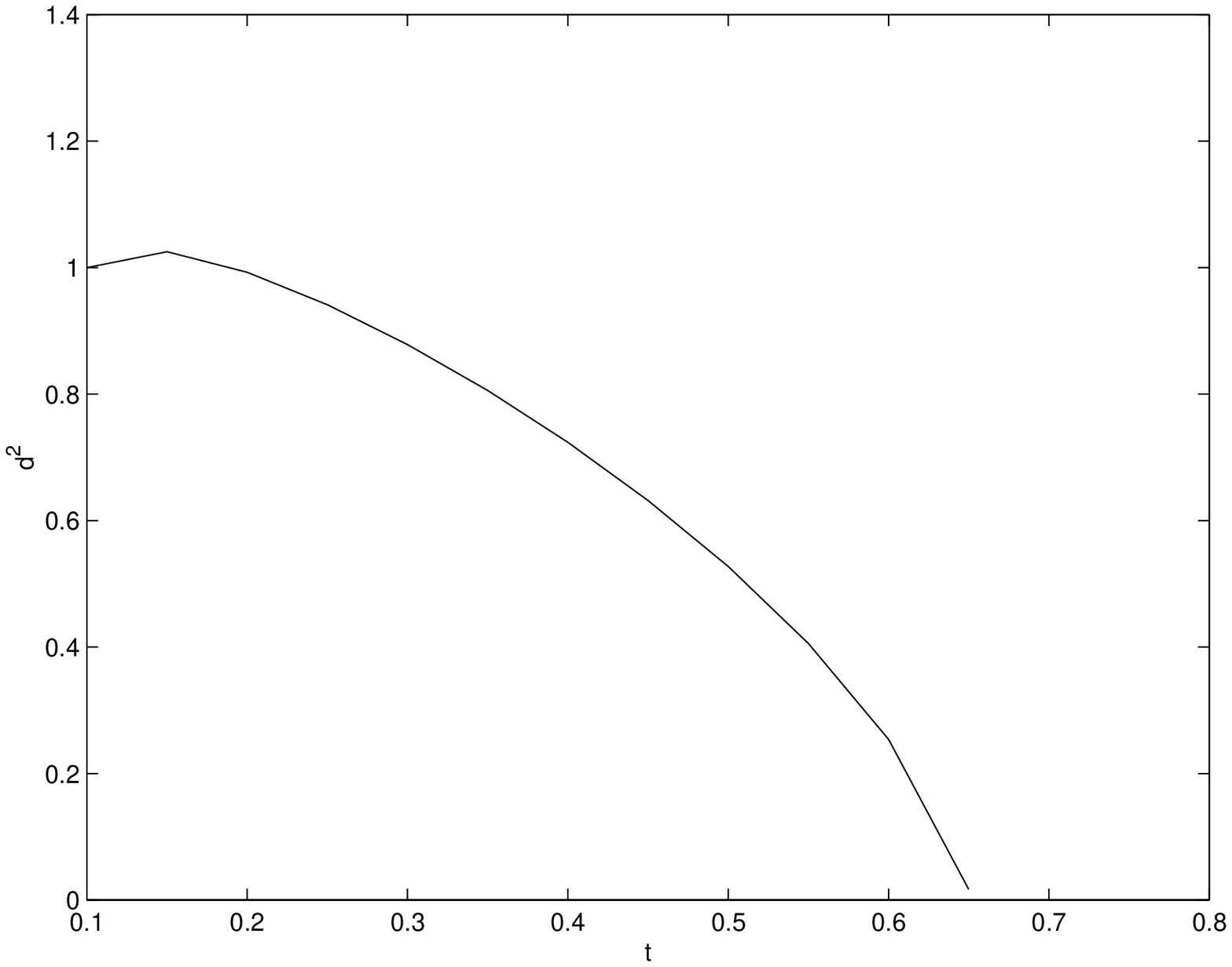}}
\centerline{Fig. 4}\par \centerline{
\includegraphics[width=3in]{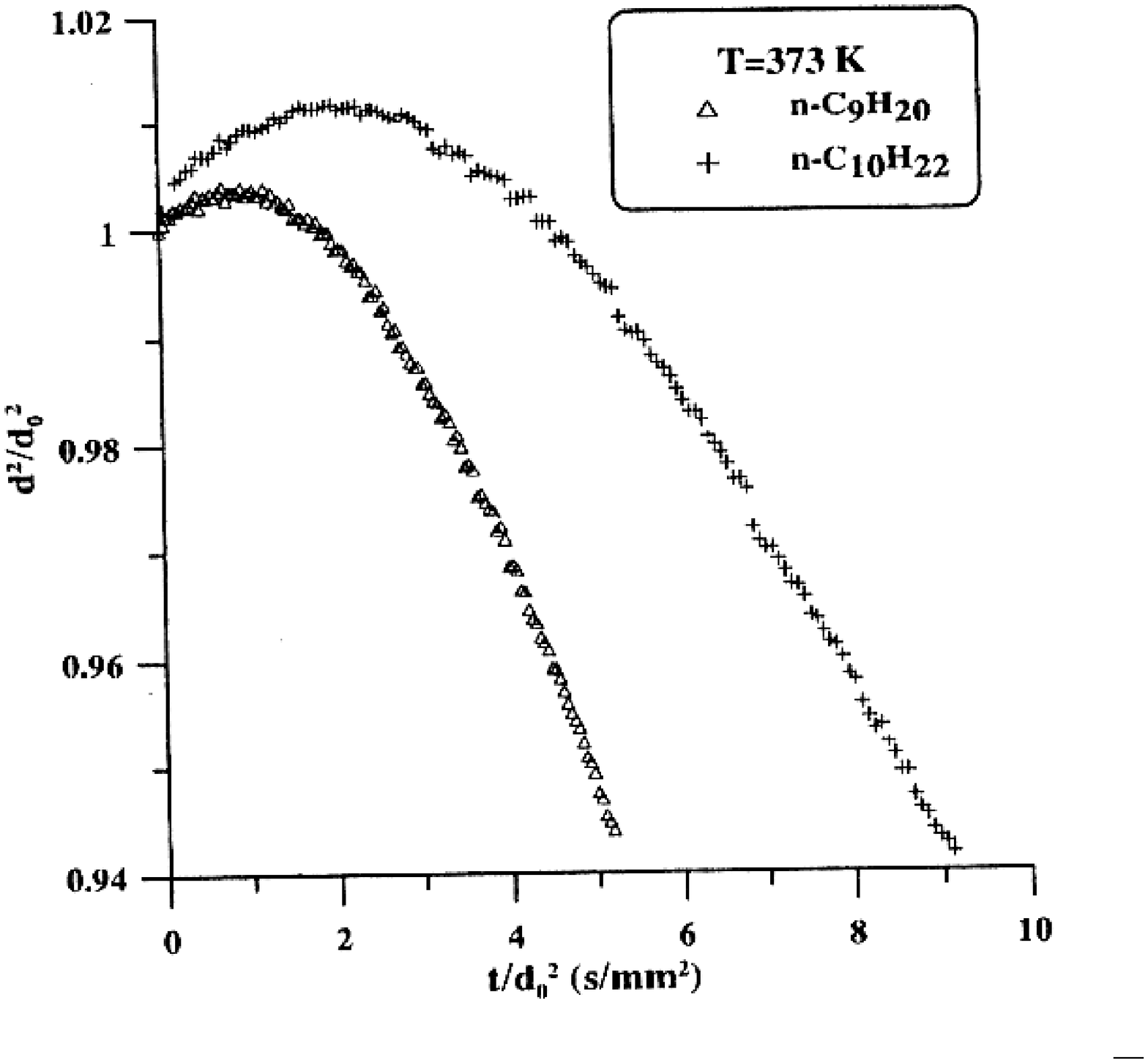}} 
\centerline{Fig. 5}\par
\smallskip
\noindent Since $ W'_{0}(\xi  )<0$ the maximal existence interval is
finite (Proposition \ref{prop1}) as can be seen in our graphic. Let
us remark that looking on the experimental curves \cite{morin} made
by the LCSR (figure 5) at the beginning, the function $
t\longrightarrow  R(t)$ is increasing around the vicinity of $ t=0.$
This fact is confirmed by our model which represents a good
improvement of our previous model\cite{goutte-1} in which the
velocity $ v_{G}(t)$ was a given function of $ t$.

\end{document}